\documentclass[12pt]{article}
\usepackage{arydshln}
\usepackage{amsmath}
\usepackage{multirow}
\usepackage{amssymb}
\usepackage{amsfonts}
\usepackage{setspace}
\usepackage{enumerate}
\usepackage{amsthm}
\usepackage{epsf}
\usepackage{graphicx}
\usepackage{changepage}   
\usepackage{cases}
\usepackage[justification=centering]{caption}
\usepackage{cite}
\def\mineappendix{
        \setcounter{section}{1}
        \setcounter{subsection}{0}
        \def\thesection{\Alph{section}}
        \def\sectionap{\@startsection  {section}{1}{\z@}
                        {-3.5ex plus-1ex minus-.2ex} {0ex plus.2ex}
                        {\reset@font\Large\bf  Appendix:  \, }
                        }
        }
\makeatother
\def\Proclaim #1. #2\par{\bigbreak\noindent{\sc#1.\enspace}{\it#2}\par}

\font\Bbbfont=msbm10
\newfam\msbfam
\textfont\msbfam=\Bbbfont \textfont\msbfam=\Bbbfont

\def\Arph{\varphi}

\def\nd{\noindent}
\newtheorem{Definition}{Definition}
\newtheorem{Theorem}{Theorem}
\newtheorem{Lemma}{Lemma}

\newtheorem{Proposition}{Proposition}
\newtheorem{Remark}{Remark}

\title{ The Almost Complex Structure on $\mathbb S^6$ and Related Schr\"odinger Flows}
\author{Qing Ding\footnote{Email: qding@fudan.edu.cn; Fax: 0086-21-65646073}
	\ \ \  and  \ \ \ Shiping Zhong \footnote{New address: School of Mathematics and Computer Sciences, Gannan Normal University, Ganzhou 431000,
P.R. China; Email: spzhong15@fudan.edu.cn}\\
	 School of Mathematical Sciences,  \\
 Fudan University, Shanghai 200433, P.R. China}
	\date{}
	
\begin{document}	
	
\maketitle

\begin{abstract}
In this paper, by using the $G_2$-structure on Im$(\mathbb O)\cong\mathbb R^7$ from the octonions $\mathbb O$, the $G_2$-binormal motion of curves $\gamma(t,s)$ in $\mathbb R^7$ associated to the almost complex structure on $\mathbb S^6$ is studied. The motion is proved to be equivalent to Schr\"odinger flows from $\mathbb R^1$ to $\mathbb S^6$, and also to a nonlinear Schr\"odinger-type system in three unknown complex functions that generalizes the famous correspondence between the binormal motion of curves in $\mathbb{R}^3$ and the focusing nonlinear Schr\"odinger equation. Some related geometric properties of the surface $\Sigma$ in Im$(\mathbb O)$ swept by $\gamma(t,s)$ are determined.
\end{abstract}


\nd Keywords: Almost complex structure, Schr\"odinger flow, $G_2$-structure, octonions

\section *{\S 1. Introduction}
The study of moving curves in a Riemannian or pseudo-Riemannian manifold, especially in the Euclidean or Minkowski spaces,
is an attractive topic in differential geometry, as it has applications in physics, such as the deformation of a thin vortex filament in inviscid fluid
\cite{Arms, Meleshko}, kinematics of interfaces in crystal growth \cite{Brower1, LangerJS}, viscous fingering in a Hele-shaw cell \cite{Saffman}, etc. It is well-known that in considering a motion of a vortex filament in $\mathbb R^3$, one encounters the Da Rios equation as follows (see \cite{Rios},\cite{Arms}):
\begin{eqnarray}
\gamma_t=\gamma_s\times \gamma_{ss},\label{1}
\end{eqnarray}
where $\gamma=\gamma(t,s)\in \mathbb R^3$ is the centerline (curve) of the vortex filament represented by a vector-valued function with respect to arclength $s$ and time $t$, the subscript stands for partial derivative with respect to the indicated variable and $\times$ denotes the cross product between vectors in $\mathbb R^3$. Eq.(\ref{1}) is also called the binormal motion of space curves in $\mathbb R^3$, since
the right-hand-side of (\ref{1}) is equal to a multiple of the curvature and the binormal vector at $\gamma$. The distinguishing features of Eq.(\ref{1}) are its completely integrability and the fact that it is equivalent to the focusing nonlinear Schr\"odinger equation (NLS):
$$
\varphi_t+\varphi_{ss}+2|\varphi|^2\varphi=0.$$
This produces the so-called Da Rios-NLS correspondence (refer to \cite{Hasimoto}, \cite{Langer}), that is, if a curve $\gamma(t,s)$ evolves according to the Da Rios equation (\ref{1}), then the associated complex function $\varphi=\kappa(t,s) \exp\left(i\int^s\tau(t,x)dx\right)$ by Hasimoto transform evolves according to the NLS equation, where $\kappa$ and $\tau$ stand for the curvature and the torsion curvature at $\gamma(t,s)$ respectively.

On the other hand, Euclidean 3-space $\mathbb R^3$ can be regarded as the imaginary part Im($\mathbb H$) of the quaternions $\mathbb H$ and the cross product $\times$ on $\mathbb R^3\cong \hbox{Im}(\mathbb H)$ is created by the quaternion algebraic structure in Im($\mathbb H$). It is also well-known that the quaternions are contained in the octonions $\mathbb O$ and thus the imaginary part Im($\mathbb H$) in the imaginary part Im($\mathbb O$) of the octonions $\mathbb O$. It is not a surprise that there is a cross product on Im($\mathbb O$)$\cong\mathbb R^7$ induced by the octonion algebraic structure, which includes the cross product on Im($\mathbb H$)$\cong\mathbb R^3$ as a special case. Therefore, as a natural generalization of Eq.(\ref{1}) in higher dimensions, the following equation in Im($\mathbb O$)$\cong\mathbb R^7$:
\begin{eqnarray}
\gamma_t=\gamma_s\times \gamma_{ss}, ~~\gamma(t,s)\in \mathbb R^7\label{2}
\end{eqnarray}
looks very interesting from the purely mathematical point of view. Eq.(\ref{2}) is in fact the $G_2$-binormal motion of curves in $\mathbb{R}^7$, as we shall see below. However, to our disappointment, we haven't found physical applications of Eq.(\ref{2}) yet in the literature.

The aim of this paper is to give a geometric interpretation of Eq.(\ref{2}) associated to the almost complex structure on $\mathbb S^6$. To our surprise, Eq.(\ref{2}) is proved to be equivalent to Schr\"odinger flows of maps from $\mathbb R^1$ to the 6-sphere $\mathbb S^6\hookrightarrow \mathbb R^7$, where $\mathbb S^6$ is equipped with the standard almost complex structure called Kirchhoff's almost complex structure that is not integrable. It was proved by Borel and Serre in 1953 in \cite{Borel} that $\mathbb S^{2n}$ admits an almost complex structure if and only if $ n = 1 $ or $3$. In 1993, Calabi and Gluck \cite{Gluck} proved that the best almost complex structure on $\mathbb S^6$ is the one constructed by Kirchhoff in the sense that it has the smallest volume in a class of sections of the bundle $O(8)=U(4)$ over $\mathbb S^6$. From the viewpoint of Schr\"odinger flows, one may obtain different equations (\ref{2}) by choosing different almost complex structures on $\mathbb S^6$ in advance. This indicates that Eq.(\ref{2}) not only is a higher dimensional generalization of Eq.(\ref{1}), but also relates to almost complex structures on $\mathbb S^{6}$. An old problem in this aspect is whether or not there is a complex structure on $\mathbb S^6$ (refer to \cite{Bryant11, LeBrun}). This gives us further motivations for studying Eq.(\ref{2}). The aim of this paper
can also be regarded as a contribution to our understanding of almost complex structures on $\mathbb S^6 $ and the $G_2$-structure on Im$(\mathbb O) = \mathbb{R}^7 $ via Schr\"odinger flows. Furthermore, Eq.(\ref{2}) is also shown to be equivalent to a nonlinear Schr\"odinger-type system in three unknown complex functions, which sets the famous Da Rios-NLS correspondence as a special case. Some geometric properties of the surface $\Sigma$ in Im$(\mathbb O)$ swept by $\gamma(t,s)$ are characterized by the unknown functions in the nonlinear Schr\"odinger-type system.

The paper is organized as follows. Section \S2 gives preliminaries about Schr\"odinger flows from a Riemannian manifold to an almost Hermitian manifold, the Cayley-Dickson construction, the exceptional simple Lie group $G_2$ and $G_2$-frame in Im($\mathbb O$)$\cong\mathbb R^7$.  In \S3, we construct the complexified $G_2$-frame and establish a related Frenet formula along curves.
In \S4, we give a proof of the correspondence between Eq.(\ref{2}) and a nonlinear Schr\"odinger-type system.
In \S5 we exploit the geometric properties of the surface swept by $G_2$-binormal moving curves in terms of the unknown-functions in the nonlinear Schr\"odinger-type system.

\section * {\S 2. Preliminaries}
In this section, we recall briefly the geometric concept of Schr\"odinger flows from a Riemannian manifold to an almost Hermitian manifold. We also give some facts about the Cayley-Dickson construction, the octonions, the exceptional simple Lie group $G_2$ and $G_2$-frame.

\subsection * {\S 2.1 Schr\"odinger flows to almost Hermitian manifolds}

The motivation for introducing Schr\"odinger flows comes from the fact that the Heisenberg model in condensed matter physics is described as the equation of Schr\"odinger flows from $\mathbb R^1$ to the 2-sphere $\mathbb S^2\hookrightarrow \mathbb R^3$,
in which the standard complex structure on $\mathbb S^2$ is used. First of all, we recall the definition of Schr\"odinger flows from the Riemannian manifold $(M,g)$ to the almost Hermitian manifold $(N,J,h)$, where $J$ is an almost complex structure compatible to the metric $h$ on $N$ (refer to \cite{dq,dw,Nahmod,TeUh}, for example). Some other related geometric flows, such as KdV geometric flows, are described in \cite{dinghe,DW1}.

\begin{Definition} A map $u=u(t,x): [0,T)\times M\to N$, where $0<
T\le\infty$, is called a Schr\"odinger flow from
$(M,g)$ to $(N,J,h)$ if $u$ satisfies the following equation of the
Hamiltonian gradient flow
\begin{eqnarray*} u_t=J_u\nabla
E(u), 
\end{eqnarray*}where $E(u)$ is the
energy functional of $u: M\to N$.
\end{Definition}
Recall that the energy $E(u)$ of $u: M\rightarrow N $ is defined by
\begin{eqnarray*}
	E(u)=\int_M e(u) dv_g,
\end{eqnarray*}
where, in a local chart $(x_{\alpha})$ of $M$, $e(u)=\frac{1}{2}g^{\alpha\beta}h_{jk}(u)\frac{\partial u^j}{\partial x_\alpha}\frac{\partial u^k}{\partial x_\beta}$.
It is easy to verify that the gradient $\nabla E(u)$ is exactly the tension field $\tau(u)$ of map $u$, which is expressed in local coordinates as
\begin{eqnarray*}
	\tau(u)^i=\Delta_Mu^i+g^{\alpha\beta}\Gamma^i_{jk}(u)\frac{\partial u^j}{\partial x_\alpha}\frac{\partial u^k}{\partial x_\beta},
\end{eqnarray*}
where $\Delta_M$ is the Laplace-Beltrami operator on $M$ and $\Gamma^i_{jk}$ are the Christoffel symbols of N. Hence, the equation of the Schr\"odinger flows of $u: M\rightarrow N$ can  be rewritten as
\begin{eqnarray}\label{flow}
	u_t=J_u\tau(u).
\end{eqnarray}

Next, it is easy to see that when the target $N =\mathbb{C}$, the complex plane, Eq.(\ref{flow}) is nothing but the linear Schrodinger equation. When $N=\mathbb S^n\hookrightarrow \mathbb R^{n+1}$ $(n=2,6)$, the tension field of map $u: M\to\mathbb S^n\hookrightarrow \mathbb R^{n+1}$ is given by $\tau(u)=\Delta u+|\nabla u|^2u$. Noting that the standard almost complex structure $J$ at $u\in \mathbb S^n$ $(n=2,6)$
is given by  $J_u=u\times : T_u \mathbb S^n \rightarrow T_u \mathbb S^n$, in which $J$ is integral only when $n=2$, we obtain the equation of  Schr\"odinger flows from $M$ to $\mathbb S^n\hookrightarrow\mathbb R^{n+1}$ as follows
\begin{eqnarray*}
	u_t=u\times\Delta u,~~ u: M\to \mathbb S^n\hookrightarrow\mathbb R^{n+1}~(n=2,6).
\end{eqnarray*}

\subsection * {\S 2.2 The Cayley-Dickson construction and $G_2$-frame }

Let $\mathbb{A}$ be an algebra over the field $\mathbb{R}$ which is not necessarily associative but finite-dimensional. A linear mapping $a \rightarrow \overline{a}$ of $\mathbb{A}$ to itself is said to be a conjugation or involutory anti-automorphism if $\overline{\overline{a}}=a$ and $\overline{ab}=\overline{b}\overline{a}$ for any elements $a,\, b\in \mathbb{A}$ (an case $\overline{a}=a$ is not excluded).

\begin{Definition}	\text{(Cayley-Dickson construction {\cite{Baez}\,,\cite{Harvey}})}
Consider the vector space of the direct sum of two copies of $\mathbb{A}$:
$\mathbb{A}^2=\mathbb{A}\oplus\mathbb{A}.$ A multiplication on $\mathbb{A}^2$ is defined as:
\begin{eqnarray*} 
(a,b)(c,d)=(ac-d\overline{b}\,,\,\overline{a}d+cb).
\end{eqnarray*}
\end{Definition}

It is easy to check that relative to the multiplication the vector space
$\mathbb{A}^2$ is an algebra of dimension $2\cdot dim(\mathbb{A})$. This is called the doubling of $\mathbb{A}$.

\begin{Remark}
The correspondence $a \rightarrow (a,0)$ is a monomorphism of $\mathbb{A}$ into $\mathbb{A}^2$. Therefore we will identify elements $a$ and $(a, 0)$ and consequently assume $\mathbb{A}$ is a subalgebra of $\mathbb{A}^2$. If $\mathbb{A}$ has an identity element, then the element $1 = (1, 0)$ is obviously an identity element in $\mathbb{A}^2$.
\end{Remark}

An important element in $\mathbb{A}^2$ is $e=(0,1)$. It follows from the definition of multiplication that $be=(0,b)$ and hence $(a,b)=a+be$ for all $a,b\in \mathbb{A}$. Thus every element of the algebra $\mathbb{A}^2$ is uniquely written as $a+be$. Moreover, as can be easily checked, the following identities are true:
\begin{eqnarray} \label{indentities}
a(be)=(ba)e,\hspace{4mm} (ae)b=(a\overline{b})e,\hspace{4mm} (ae)(be)=-\overline{b}a.
\end{eqnarray}
In particular $e^2=-1$.

To iterate the Cayley-Dickson construction it is necessary to define a conjugation in $\mathbb{A}^2$. This will be done by the formula
\begin{eqnarray*}
	\overline{a+b}=\overline{a}-be.
\end{eqnarray*}
The doubling $\mathbb{R}^2$ of the field $\mathbb{R}$ is the algebra $\mathbb{C}$ of complex numbers and the doubling $\mathbb{C}^2$ of $\mathbb{C}$ is the algebra of quaternions $\mathbb{H}$. In the latter case $e$ is denoted by $j$ and $ie$ is denoted by $k$, so a general quaternion is of the form $r = r_1 +r_2i+r_3j +r_4k$, where $r_i\in \mathbb{R}, i=1, 2, 3, 4$. Due to the identity (\ref{indentities}) $ea=\overline{a}e$ for all $a \in \mathbb{A}$, one may verify that $\mathbb{H}$ is not commutative.

The doubling algebra $\mathbb O=\mathbb{H}^2$ of $\mathbb H$ is the Cayley algebra which is not commutative and associative, and its elements are called octonions or Cayley numbers.
By definition every octonion is of the form $\xi=a+be$, where $a$ and $b$ are quaternions. The basis of $\mathbb O$ consists of $\{1, 	i,\,j,\,k,\,l,\,il,\,jl,\,kl\}$ in which we replace $e$ by $l$. The square of each of these elements is $-1$ except the unit element $1$. The full multiplication table is summarized in Table \ref{o123}.
\begin{table}[!h]
\centering
	\begin{tabular}{|c|c|c|c|c|c|c|c|}
		\hline
		$$ & $i$ & $j$ & $k$  & $l$ & $il$ & $jl$ & $kl$\\
		\hline
		$i$ & $-1$ & $k$ & $-j$ &
		$il$ & $-l$ & $-kl$ & $jl$ \\
		\hline
		$j$ & $-k$ & $-1$ & $i$&
		$jl$ & $kl$ & $-l$ & $-il$ \\
		\hline
		$k$ & $j$ & $-i$ & $-1$&
		$kl$ & $-jl$ & $il$ & $-l$ \\
		\hline
		$l$ & $-il$ & $-jl$ & $-kl$&
		$-1$ & $i$ & $j$ & $k$ \\
		\hline
		$il$ & $l$ & $-kl$ & $jl$&
		$-i$ & $-1$ & $-k$ & $j$ \\
		\hline
		$jl$ & $kl$ & $l$ & $-il$&
		$-j$ & $k$ & $-1$ & $-i$ \\
		\hline
		$kl$ & $-jl$ & $il$ & $l$&
		$-k$ & $-j$ & $i$ & $-1$ \\
		\hline	
        \end{tabular}
        \caption{The multiplication table of $\mathbb O$} \label{o123}
\end{table}	

The  group $G_2$ is defined to be the automorphism group of the octonions $\mathbb O :$
\begin{eqnarray*}
G_2=\{g\in Iso_{\mathbb{R}}\mathbb O\,\,|\,\,g(xy)=g(x)g(y),\,\, \forall \,\,x, y\in \mathbb O\},
\end{eqnarray*}
where Iso$_{\mathbb{R}}\mathbb O=O(8)$ denotes the set of all $\mathbb{R}$-linear isomorphisms of $\mathbb O$.
The cross product $x\times y$ and scalar product $\langle x,y\rangle$¡± of $\mathbb O$ are determined respectively by the multiplicity on $\mathbb O$
as follows:
\begin{eqnarray*}
x\times y=(1/2)(\bar{y}x-\bar{x}y),\,\,\langle x,y\rangle=(1/2)(\bar{x}y+\bar{y}x),
\end{eqnarray*}
where $\bar{x}=2\langle x,1\rangle-x$ is the conjugation of $x\in\mathbb O$. One may verify directly that $x\times y\in $ Im$(\mathbb O)$, $\forall x,\, y\in\, $Im$(\mathbb O)$,\,where Im$(\mathbb O)=\{x\in \mathbb O\,|\,\langle x,1\rangle=0\}$. This induces a cross product $\times$ among vectors in Im$(\mathbb O)\cong \mathbb R^7$.

The differential version of $G_2$ realization leads to the so-called $G_2$-frame of $ \mathbb{R}^7 = $ Im$(\mathbb O)$. In 1958, Calabi\cite{Calabi} first constructed $G_2$-structure equations of a submanifold in $\mathbb{R}^7 $. In 1982, Bryant \cite{Bryant,Bryant11} gave a more concrete representation of $G_2$ by taking account of the algebraic properties of the octonions $ \mathbb O$.
	
From the standard basis of Im$(\mathbb O)=\text{Span}_{\mathbb{R}}\{	i,\,j,\,k,\,l,\,il,\,jl,\,kl\}$,
we define a basis of the complexification of Im$(\mathbb O)$ over $ \mathbb{C} $:
\begin{eqnarray*}\label{G2basis1}
\begin{split}
 N=\frac{1}{\sqrt{2}}(1-\sqrt{-1}l),\,\,\overline{N}=\frac{1}{\sqrt{2}}(1+\sqrt{-1}l),\hspace{2cm}\\ E_1=iN,\,\,E_2=jN,\,E_3=-kN,\,\,\overline{E}_1=i\,\overline{N},\,\,\overline{E}_2=j\,\,\overline{N},\,\overline{E}_3=-k\,\,\overline{N}.
\end{split}
\end{eqnarray*}
A basis $(e_4\,\,\,\,f\,\,\,\,\overline{f})$ of $\mathbb{C}\otimes_{\mathbb{R}} \mathbb O$ is said to be admissible, if there exists $g\in G_2$ such that $(e_4\,\,\,\,f\,\,\,\,\overline{f})^T=g(l\,\,\,\,E\,\,\,\,\overline{E})^T$, where $E=(E_1,E_2,E_3)$.
Usually, $(e_4\,\,\,\,f\,\,\,\,\overline{f})$ is called a complexified $G_2$-frame.
\begin{Theorem}\text{(Bryant\cite{Bryant})}\label{th-ab}
For a complexified $G_2$-frame $(e_4\,\,\,\,f\,\,\,\,\overline{f})$, we have
\begin{eqnarray*}
\begin{array}{c@{\hspace{-5pt}}l}
d\left(\begin{array}{ccc;{2pt/2pt}ccc}
e_4\\
f\\
\overline{f}
\end{array}\right)
=\left(
\begin{array}{@{}c|c|c@{}}
0 & -\sqrt{2}\sqrt{-1}\,\,\theta & \sqrt{2}\sqrt{-1}\,\,\overline{\theta} \\
\hline
-\sqrt{2}\sqrt{-1}\,\,\overline{\theta}^T & \kappa & [\theta] \\
\hline
\sqrt{2}\sqrt{-1}\,\,\theta^T & \overline{[\theta]} & \overline{\kappa}
\end{array}
\right)
\end{array}
\left(\begin{array}{ccc;{2pt/2pt}ccc}
e_4\\
f\\
\overline{f}
\end{array}\right),
\end{eqnarray*}
where $\theta=(\theta^1\,\, \theta^2\,\,\theta^3)$ is an $M_{1\times3}(\mathbb{C})$ valued 1-form, $\kappa$ is an $\textit{su}(3)$
valued 1-form which satisfies
\begin{eqnarray*}
\kappa+\bar{\kappa}^T=0_{3\times3},\,\,\,\, tr \kappa=0
\end{eqnarray*}%
and
\begin{eqnarray*}
\begin{split}
[\theta]=
\left(
\begin{array}{@{}ccc@{}}
0 & -\theta^3 & \theta^2 \\
\theta^3 & 0 & -\theta^1 \\
-\theta^2 & \theta^1 & 0
\end{array}
\right).
\end{split}
\end{eqnarray*}%
\end{Theorem}

The above structure of the $G_2$-frame will play a crucial role in the proof of Theorem 4, as we shall see below.
We should mention that the standard Kirchhoff's complex structure $J$ on $\mathbb S^6 $ is explicitly given as follows: for $u\in \mathbb S^6$,		
\begin{eqnarray*}
J_u: T_u \mathbb S^6 \rightarrow T_u \mathbb S^6, \,\,X\mapsto J_u(X)=u\times X,\,\,\,\,X\in \mathbb T_u\mathbb S^6.
\end{eqnarray*}

\section * {\S 3 $G_2$-structure equations of curves in Im($\mathbb O$)}
In this section, we describe the construction of $G_2$-frame along curves in Im$(\mathbb O)\cong \mathbb R^7$ and then take its complexification. One may refer to  \cite{Hashimoto,ohashi2011,ohashi20132} for details. Based on the complexification, we present a Frenet formula
of the complexified $G_2$-frame along a curve in Im$(\mathbb O)\cong \mathbb R^7$.

Let $\gamma(s)$ be a unit speed curve in Im$(\mathbb O)\cong \mathbb R^7$. We set $k_1(s)=\parallel\gamma_{ss}(s)\parallel$ and assume that this function does not vanish anywhere. Now we define a $G_2$-frame along the curve as follows
\begin{eqnarray*}\label{G200}
\begin{split}
	I_4(s)&=\gamma_s(s),\,\,\,I_1(s)=\frac{1}{k_1}I_{4s},\,\,\,I_5(s)=I_1\times I_4,\\
	I_2(s)&=\frac{1}{\kappa_2}(I_{1s}-\langle I_{1s},I_4\rangle I_4-\langle I_{1s},I_5\rangle I_5),\\
	I_3(s)&=I_1\times I_2,\,\,\,I_6(s)=I_2\times I_4,\,\,\,I_7(s)=I_3\times I_4,
\end{split}
\end{eqnarray*}
where
\begin{eqnarray*}\label{G201}
	\kappa_2(s)=\sqrt{\parallel I_{1s}\parallel^2-\langle I_{1s},I_4\rangle ^2-\langle I_{1s},I_5\rangle ^2}> 0
\end{eqnarray*}
is assumed and $I_5=I_1\times I_4$ is usually regarded as the $G_2$-binormal vector along the curve.

The multiplication table of $(I_4 \,\, I_1 \,\, I_2 \,\, I_3 \,\, I_5 \,\, I_6 \,\, I_7)$ coincides with that of $(l \,\, i \,\, j \,\, k \,\, il \,\, jl \,\, kl)$. In other words, there exists a $G_2$-valued function $g$ such that
\begin{eqnarray*}
	(I_4 \,\, I_1 \,\, I_2 \,\, I_3 \,\, I_5 \,\, I_6 \,\, I_7)=(g(l) \,\, g(i) \,\, g(j) \,\, g(k) \,\, g(il) \,\, g(jl) \,\, g(kl)).
\end{eqnarray*}

If $\kappa_2=0$ and we take $I_2(s)\in $ (span$_\mathbb{R} \{I_4, I_1, I_5\})^\bot$ with $|I_2(s)|=1$,
$I_3(s)=I_1\times I_2,\,\,\,I_6(s)=I_2\times I_4$, and $I_7(s)=I_3\times I_4,$ then $(I_4 \,\, I_1 \,\, I_2 \,\, I_3 \,\, I_5 \,\, I_6 \,\, I_7)$ also
consists of a $G_2$-frame along the curve in Im$(\mathbb O)\cong \mathbb R^7$, in which $\{I_4,\,I_1,\,I_5\}$ consists of an autonomy system,
that is, it satisfies the formula
\begin{eqnarray*}
\begin{array}{c@{\hspace{-5pt}}l}
\left(\begin{array}{ccc}
I_4\\
I_1\\
I_5
\end{array}\right)_s
=\left(
\begin{array}{@{}ccc@{}}
0 & k_1 & 0 \\
-k_1 & 0 & \rho_1 \\
0 & -\rho_1 & 0
\end{array}
\right)
\end{array}
\left(\begin{array}{ccc;{2pt/2pt}ccc}
I_4\\
I_1\\
I_5
\end{array}\right),
\end{eqnarray*}
where $\rho_1=\langle I_{1s},I_5\rangle.$

\begin{Proposition}(\cite{ohashi2011,ohashi20132} ) \label{Pro1}
	Let $\gamma: I=(0,1) \rightarrow $ Im($\mathbb O$) be a curve with $k_1 >0$. The associated $G_2$-frame $(I_4 \,\, I_1 \,\, I_2 \,\, I_3 \,\, I_5 \,\, I_6 \,\, I_7)$ satisfies the following differential equation
	\begin{eqnarray}\label{G2basis}
	\begin{array}{c@{\hspace{-5pt}}l}
	\left(\begin{array}{ccc;{2pt/2pt}ccc}
	I_4\\
	I_1\\
	I_2\\
	I_3\\
	I_5\\
	I_6\\
	I_7
	\end{array}\right)_s
	=\left(
	\begin{array}{@{}c|ccc|ccc@{}}
	0 & k_1 & 0 & 0 & 0 & 0 & 0 \\
	\hline
	-k_1 & 0 & \kappa_2 & 0 & \rho_1 & 0& 0 \\
	0 & -\kappa_2 & 0 & \alpha & 0 & \rho_2 & \beta_1 \\
	0 & 0 & -\alpha & 0 & 0 & \beta_2 & \rho_3 \\
	\hline
	0 & -\rho_1 & 0 & 0 & 0 & \kappa_2 & 0 \\
	0 & 0 & -\rho_2 & -\beta_2 & -\kappa_2 & 0 & \alpha \\
	0 & 0 & -\beta_1 & -\rho_3 & 0 & -\alpha & 0 \\
	\end{array}
	\right)
	\end{array}
	\left(\begin{array}{ccc;{2pt/2pt}ccc}
	I_4\\
	I_1\\
	I_2\\
	I_3\\
	I_5\\
	I_6\\
	I_7
	\end{array}\right)
	\end{eqnarray}
	with $\rho_1=\langle I_{1s},I_5\rangle ,\, \rho_2=\langle I_{2s},I_6\rangle ,\, \rho_3=\langle I_{3s},I_7\rangle ,\, \alpha=\langle I_{2s},I_3\rangle ,\,
	\beta_1=\langle I_{2s},I_7\rangle$ and $\beta_2=\langle I_{3s},I_6\rangle $. These functions satisfy
	\begin{eqnarray}
	\rho_1+\rho_2+\rho_3&=&0 ,\label{guangxi1}\\
	\beta_1-\beta_2+k_1&=&0 .\label{guangxi2}
	\end{eqnarray}
\end{Proposition}

\begin{Remark} One notes from (\ref{G2basis}) that Eq.(\ref{2}) can be rewritten as $\gamma_t=-k_1I_5$. Hence, Eq.(\ref{2}) is the $G_2$-binormal motion of curves in $\mathbb{R}^7$. We mention that
for a given curve $\gamma(s)$ in Im($\mathbb O)\cong \mathbb R^7$, one may also establish the associated $SO(7)$-frame $\{V_1,\,V_2,\,\cdots,\,V_7\}$ and corresponding Frenet formula, as was done by Ohashi in \cite{ohashi2015}. For the relationship between the $G_2$-frame and $SO(7)$-frame,
Ohashi proved in \cite{ohashi2015} that
\begin{eqnarray*}
I_5=\cos\sigma\,V_3+\sum\limits_{i=1}^4a_{i1}\,V_{3+i},
\end{eqnarray*}
where
\begin{eqnarray*}
\cos \sigma=\langle V_3, V_3\times V_1\rangle,\,\,a_{i1}=\langle V_{i+3}, V_2\times V_1 \rangle,i\in\{1,2,3,4\}.
\end{eqnarray*}
Therefore, from the view point of $SO(7)$-frame, Eq.(\ref{2}) is a motion of curves along the direction in a combination of $\{V_3,V_4,V_5,V_6,V_7\}$ in $\mathbb{R}^7$.
\end{Remark}

The six functions $(k_1,\,\,\kappa_2,\,\,\rho_1,\,\,\rho_3,\,\,\alpha,\,\,\beta_1)$ are the complete $G_2$-invariants of $\gamma$.
Now, we give the complexification of the $G_2$-frame $(I_4 \,\, I_1 \,\, I_2 \,\, I_3 \,\, I_5 \,\, I_6 \,\, I_7)$  along $\gamma(s)$ according to Bryant in \cite{Bryant} as follows. Let
\begin{eqnarray}\label{basis}
\begin{split}
&&e_4&=&I_4,&&\\
e_1&=&r(I_1-\sqrt{-1}I_5)&,&\overline{e}_1&=&\overline{r}(I_1+\sqrt{-1}I_5),\\
e_2&=&q(I_2-\sqrt{-1}I_6)&,&\overline{e}_2&=&\overline{q}(I_2+\sqrt{-1}I_6),\\
e_3&=&-p(I_3-\sqrt{-1}I_7)&,&\overline{e}_3&=&-\overline{p}(I_3+\sqrt{-1}I_7),
\end{split}
\end{eqnarray}
where
\begin{eqnarray}\label{rqp}
\begin{split}
r&=\frac{1}{\sqrt{2}}\exp(-\sqrt{-1}\int_0^s\rho_1d\widetilde{s}),\\
q&=\frac{1}{\sqrt{2}}\exp(-\sqrt{-1}\int_0^s\rho_2d\widetilde{s}),\\
p&=\frac{1}{\sqrt{2}}\exp(-\sqrt{-1}\int_0^s\rho_3d\widetilde{s})=\sqrt{2}\,\overline{q}\,\overline{r}
\end{split}
\end{eqnarray}
and bar denotes the complex conjugation of elements in $\mathbb{C}\otimes_{\mathbb{R}} \mathbb O$, namely,
\begin{eqnarray*}
\overline{x_1+\sqrt{-1}x_2}=x_1-\sqrt{-1}x_2,\,\,\,\,\,\,x_1,x_2\in \mathbb O.
\end{eqnarray*}
One notes that the complex conjugation here is different from the conjugation over $\mathbb O$. The complex-conjugate on $\mathbb{C}\otimes_{\mathbb{R}} \mathbb O$ satisfies
\begin{eqnarray*}
\overline{xy}=\overline{x}\,\,\overline{y},\,\,\,\,\, x, y\in \mathbb{C}\otimes_\mathbb{R}\mathbb O.
\end{eqnarray*}
In the sequel, the bar is used to denote the complex conjugation on $\mathbb{C}\otimes_{\mathbb{R}}\mathbb O$, unless otherwise specified.

Eq.(\ref{basis}) can be rewritten as
\begin{eqnarray}\label{basis1}
\begin{array}{c@{\hspace{-5pt}}l}
\left(\begin{array}{ccc}
 e_4\\
  f\\
 \overline{f}
\end{array}\right)
=\left(
\begin{array}{@{}ccc@{}}
1 & 0 & 0 \\
0 & A & -\sqrt{-1}A \\
0 & \overline{A} & \sqrt{-1}\,\,\overline{A}
\end{array}
\right)
\end{array}
\left(\begin{array}{ccc;{2pt/2pt}ccc}
 I_4\\
 J_1\\
 J_2
\end{array}\right):=N_1
\left(\begin{array}{ccc;{2pt/2pt}ccc}
 I_4\\
 J_1\\
 J_2
\end{array}\right),
\end{eqnarray}
where
\begin{eqnarray}\label{fenliang1}
f=\left(\begin{array}{ccc;{2pt/2pt}ccc}
 e_1\\
 e_2\\
 e_3
\end{array}\right),\,\,\,\,
J_1=\left(\begin{array}{ccc;{2pt/2pt}ccc}
 I_1\\
 I_2\\
 I_3
\end{array}\right),\,\,\,\,
J_2=\left(\begin{array}{ccc;{2pt/2pt}ccc}
 I_5\\
 I_6\\
 I_7
\end{array}\right),\,\,\,\,
A=\left(\begin{array}{@{}ccc@{}}
r & 0 & 0 \\
0 & q & 0 \\
0 & 0 & -p
\end{array}\right).
\end{eqnarray}

\begin{Theorem}\label{thembasis}
For the complexified $G_2$-frame $(e_4 \,\, e_1 \,\, e_2 \,\, e_3 \,\, \overline{e}_1 \,\, \overline{e}_2 \,\, \overline{e}_3)$ along $\gamma$, we have the following Frenet formula:
\begin{eqnarray}\label{fuhaibasis}
\begin{array}{c@{\hspace{-5pt}}l}
\left(\begin{array}{ccc;{2pt/2pt}ccc}
 e_4\\
 e_1\\
 e_2\\
 e_3\\
 \overline{e}_1\\
 \overline{e}_2\\
 \overline{e}_3
\end{array}\right)_s
=\left(
\begin{array}{@{}c|ccc|ccc@{}}
0 & \Arph_1 & 0 & 0 & \overline{\Arph}_1 & 0 & 0 \\
\hline
-\overline{\Arph}_1 & 0 & \Arph_2 & 0 & 0 & 0& 0 \\
0 & -\overline{\Arph}_2 & 0 & \Arph_3 & 0 & 0 & -\frac{\sqrt{-1}}{\sqrt{2}}\Arph_1 \\
0 & 0 & -\overline{\Arph}_3 & 0 & 0 & \frac{\sqrt{-1}}{\sqrt{2}}\Arph_1 & 0 \\
\hline
-\Arph_1 & 0 & 0 & 0 & 0 & \overline{\Arph}_2 & 0 \\
0 & 0 & 0 & \frac{\sqrt{-1}}{\sqrt{2}}\overline{\Arph}_1 & -\Arph_2 & 0 & \overline{\Arph}_3 \\
0 & 0 & -\frac{\sqrt{-1}}{\sqrt{2}}\overline{\Arph}_1 & 0 & 0 & -\Arph_3 & 0 \\
\end{array}
\right)
\end{array}
\left(\begin{array}{ccc;{2pt/2pt}ccc}
 e_4\\
 e_1\\
 e_2\\
 e_3\\
 \overline{e}_1\\
 \overline{e}_2\\
 \overline{e}_3
\end{array}\right),
\end{eqnarray}
where
\begin{eqnarray}\label{Arph1}
\Arph_1=k_1\overline{r},\,\,\,\,\,\Arph_2=2\kappa_2r\overline{q},\,\,\,\,\,\Arph_3=-\sqrt{2}q^2r[2\alpha+\sqrt{-1}(\beta_1+\beta_2)].
\end{eqnarray}
\end{Theorem}

{\bf Proof}:
First of all, we can rewrite (\ref{G2basis}) as
\begin{eqnarray*}\label{basis3}
\begin{array}{c@{\hspace{-5pt}}l}
\left(\begin{array}{ccc}
 I_4\\
 J_1\\
 J_2
\end{array}\right)_s
=\left(
\begin{array}{@{}ccc@{}}
0 & u & 0 \\
-u^T & B & C \\
0 & -C^T & B
\end{array}
\right)
\end{array}
\left(\begin{array}{ccc;{2pt/2pt}ccc}
 I_4\\
 J_1\\
 J_2
\end{array}\right):=N_2
\left(\begin{array}{ccc;{2pt/2pt}ccc}
 I_4\\
 J_1\\
 J_2
\end{array}\right),
\end{eqnarray*}
where
\begin{eqnarray*}\label{fenliang2}
u=\left(\begin{array}{cc}
 k_1\\
 0\\
 0
\end{array}\right)^T,\,\,\,\,
B=\left(\begin{array}{cccc}
0 & \kappa_2 & 0 \\
-\kappa_2 & 0 & \alpha \\
0 & -\alpha & 0
\end{array}\right),\,\,\,\,
C=\left(\begin{array}{@{}ccc@{}}
\rho_1 & 0 & 0 \\
0 & \rho_2 & \beta_1 \\
0 & \beta_2 & \rho_3
\end{array}\right).
\end{eqnarray*}
From (\ref{basis1}) and (\ref{basis3}), we obtain
\begin{eqnarray*}
(e_4 \,\,\,f\,\,\,\overline{f}\,\,)^T_s
=(N_{1s}+N_1N_2)N_1^{-1}(e_4 \,\,\,f\,\,\,\overline{f}\,\,)^T,
\end{eqnarray*}
with
\begin{eqnarray*}
\begin{split}
&(N_{1s}+N_1N_2)N_1^{-1}
=\left(\begin{array}{cccc}
0 & u\overline{A} & uA \\
-Au^T & \eta_1 & \eta_2 \\
-\overline{A}u^T & \overline{\eta}_2 & \overline{\eta}_1
\end{array}\right),\\
\eta_1&=[2A_s+2AB+\sqrt{-1}A(C^T+C)]\overline{A},\\
\eta_2&=-\sqrt{-1}A(C-C^T)A.
\end{split}
\end{eqnarray*}
It follows from (\ref{guangxi1}),\,\,(\ref{guangxi2}),\,\,(\ref{fenliang1}) and  (\ref{fenliang2}) that
\begin{eqnarray*}
\begin{split}
u\overline{A}&=(k_1\overline{r}\,\,\,\,\,\,0\,\,\,\,\,\,0),\\
\eta_1
&=\left(\begin{array}{cccc}
0 & 2\kappa_2r\overline{q} & 0 \\
-2\kappa_2q\overline{r} & 0 & -\sqrt{2}q^2r[2\alpha+\sqrt{-1}(\beta_1+\beta_2)] \\
0 & \sqrt{2}\,\overline{q}^2\,\overline{r}[2\alpha-\sqrt{-1}(\beta_1+\beta_2)] & 0
\end{array}\right),\\
\eta_2&=\left(\begin{array}{cccc}
0 & 0 & 0 \\
0 & 0 & -\frac{\sqrt{-1}}{\sqrt{2}}k_1\overline{r} \\
0 & \frac{\sqrt{-1}}{\sqrt{2}}k_1\overline{r} & 0
\end{array}\right).
\end{split}
\end{eqnarray*}
By setting
\begin{eqnarray*}
\Arph_1=k_1\overline{r},\,\,\,\,\,\Arph_2=2\kappa_2r\overline{q}~~\,\hbox{and}~~\,\Arph_3=-\sqrt{2}q^2r[2\alpha+\sqrt{-1}(\beta_1+\beta_2)],
\end{eqnarray*}
we arrive at (\ref{fuhaibasis}). \qed

\medskip
Based on the complexified $G_2$-frame $(e_4 \,\, e_1 \,\, e_2 \,\, e_3 \,\, \overline{e}_1 \,\, \overline{e}_2 \,\, \overline{e}_3)$ along the curve $\gamma(s)$, we shall establish the multiplication and cross product tables with respect to $(e_4 \, e_1 \, e_2 \, e_3 \, \overline{e}_1 \, \overline{e}_2 \, \overline{e}_3)$, which are well-suited for the
development in the next section. We have the following tables, in which $e_0=-1-\sqrt{-1}e_4$:
\begin{table}[!h]
\centering
\begin{tabular}{|c|c|c|c|c|c|c|c|}
\hline
$A\backslash B$ & $e_4$ & $e_1$ & $e_2$ &  $e_3$ &  $\overline{e}_1$ & $\overline{e}_2$ &  $\overline{e}_3$\\
\hline
$e_4$ &
$ -1 $ & $-\sqrt{-1}e_1$ & $-\sqrt{-1}e_2$ &  $-\sqrt{-1}e_3$ &
$\sqrt{-1}\overline{e}_1$ & $\sqrt{-1}\overline{e}_2$ &  $\sqrt{-1}\overline{e}_3$\\
\hline
$e_1$ &
$ \sqrt{-1}e_1 $ & $0$ & $-\sqrt{2}\overline{e}_3$ &  $\sqrt{2}\overline{e}_2$ &
$e_0$ & $0$ &  $0$\\
\hline
$e_2$ &
$ \sqrt{-1}e_2 $ & $\sqrt{2}\overline{e}_3$ & $0$ &  $-\sqrt{2}\overline{e}_1$ &
$0$ & $e_0$ &  $0$\\
\hline
$e_3$ &
$ \sqrt{-1}e_3 $ & $-\sqrt{2}\overline{e}_2$ & $\sqrt{2}\overline{e}_1$ &  $0$ &
$0$ & $0$ &  $e_0$\\
\hline
$\overline{e}_1$ &
$-\sqrt{-1}\overline{e}_1 $ &$\overline{e}_0$ & $0$ &$0$ &
$0$ & $-\sqrt{2}e_3$ &  $\sqrt{2}e_2$ \\
\hline
$\overline{e}_2$ &
$ -\sqrt{-1}\overline{e}_2 $& $0$ & $\overline{e}_0$ &$0$ &
$\sqrt{2}e_3$ & $0$ &  $-\sqrt{2}e_1$ \\
\hline
$\overline{e}_3$ &
$-\sqrt{-1}\overline{e}_3 $ &$0$ & $0$ &$\overline{e}_0$&
 $-\sqrt{2}e_2$ & $\sqrt{2}e_1$ &  $0$ \\
\hline
\end{tabular}
\caption{The multiplication table $AB$}  \label{table1}
\end{table}
\begin{table}[!h]
\centering
\begin{tabular}{|c|c|c|c|c|c|c|c|}
\hline
$A\backslash B$ & $e_4$ & $e_1$ & $e_2$ &  $e_3$ &  $\overline{e}_1$ & $\overline{e}_2$ &  $\overline{e}_3$\\
\hline
$e_4$ &
$ 0 $ & $-\sqrt{-1}e_1$ & $-\sqrt{-1}e_2$ &  $-\sqrt{-1}e_3$ &
$\sqrt{-1}\overline{e}_1$ & $\sqrt{-1}\overline{e}_2$ &  $\sqrt{-1}\overline{e}_3$\\
\hline
$e_1$ &
$ \sqrt{-1}e_1 $ & $0$ & $-\sqrt{2}\overline{e}_3$ &  $\sqrt{2}\overline{e}_2$ &
$-\sqrt{-1}e_4$ & $0$ &  $0$\\
\hline
$e_2$ &
$ \sqrt{-1}e_2 $ & $\sqrt{2}\overline{e}_3$ & $0$ &  $-\sqrt{2}\overline{e}_1$ &
$0$ & $-\sqrt{-1}e_4$ &  $0$\\
\hline
$e_3$ &
$ \sqrt{-1}e_3 $ & $-\sqrt{2}\overline{e}_2$ & $\sqrt{2}\overline{e}_1$ &  $0$ &
$0$ & $0$ &  $-\sqrt{-1}e_4$\\
\hline
$\overline{e}_1$ &
$-\sqrt{-1}\overline{e}_1 $ &$\sqrt{-1}e_4$ & $0$ &$0$ &
$0$ & $-\sqrt{2}e_3$ &  $\sqrt{2}e_2$ \\
\hline
$\overline{e}_2$ &
$ -\sqrt{-1}\overline{e}_2 $& $0$ & $\sqrt{-1}e_4$ &$0$ &
$\sqrt{2}e_3$ & $0$ &  $-\sqrt{2}e_1$ \\
\hline
$\overline{e}_3$ &
$-\sqrt{-1}\overline{e}_3 $ &$0$ & $0$ &$\sqrt{-1}e_4$&
$-\sqrt{2}e_2$ & $\sqrt{2}e_1$ &  $0$ \\
\hline
\end{tabular}
\caption{The multiplication table of the cross product $A\times B$} \label{table2}
\end{table}

Furthermore, the complexified $G_2$-frame $(e_4 \,\, f \,\, \overline{f})^T$  satisfies
\begin{eqnarray}
\langle e_4,e_i \rangle=0,\,\,\, \langle e_i,e_j \rangle=\langle \overline{e}_i,\overline{e}_j \rangle=0,\,\,\, \langle e_i,\overline{e}_j \rangle=\delta_{ij},\label{neiji1}\\
e_i\times e_4=\sqrt{-1}e_i,\,\,\,\,\,\,\langle e_1\times e_2,e_3 \rangle=-\sqrt{2}\label{neiji2},
\end{eqnarray}
for any $i\in\{1,2,3\}$. By using Table \ref{table1},\,\,Table \ref{table2},\,\,(\ref{neiji1}) and (\ref{neiji2}), we may directly deduce Theorem \ref{th-ab}.

\section * {\S 4 Schr\"odinger flows to 6-sphere }

We have already seen that the equation of Schr\"odinger flows from $ \mathbb{R}^1 $ to
$ N =\mathbb S^2 \hookrightarrow \mathbb R^3$ is actually the Heisenberg ferromagnet model:
\begin{eqnarray*}
T_t=T\times T_{ss},
\end{eqnarray*}
which is equivalent to the Da Rios equation (\ref{1}).
For $\mathbb S^6=\{(x_1,\,\,x_2,\,\,x_3,\,\,x_4,\,\,x_5,\,\,x_6,\,\,x_7):\,\,\sum\limits_{i=1}^{7} x_i^2=1\}$ in $\mathbb{R}^7$, we know that
the equation (\ref{flow}) of Schr\"odinger flows from $\mathbb{R}^1$ to $(\mathbb S^6,J)$ reads
\begin{eqnarray} \label{flows6}
u_t=u\times u_{ss},
\end{eqnarray}
where $u=(u_1,\,\,u_2,\,\,u_3,\,\,u_4,\,\,u_5,\,\,u_6,\,\,u_7)\in \mathbb{R}^7$ with $\sum\limits_{i=1}^{7} u_i^2=1$.


Returning to Eq.(\ref{2}), we have
\begin{Proposition}
Suppose that $\gamma(t,s)$ evolves according to Eq.(\ref{2}). Then the arclength parameter $s$ is independent of time $t$ for all $t>0$.
\end{Proposition}

{\bf Proof}:
It suffices to prove that $\frac{d}{dt}|\gamma_s|^2=0$. In fact, from Eq.(\ref{2}), we have
\begin{eqnarray*}
\frac{d}{dt}|\gamma_s|^2=\langle \gamma_s, \gamma_s \rangle_t=2\langle \gamma_{st}, \gamma_s \rangle=2\langle \gamma_{ts}, \gamma_s \rangle=2\langle (\gamma_s\times\gamma_{ss})_s,\gamma_s\rangle=2\langle \gamma_{s}\times\gamma_{sss},\gamma_s\rangle=0. \qed
\end{eqnarray*}

Let $e_4=\gamma_s$ and we obtain from Eq.(\ref{2}) that
\begin{eqnarray*}
e_{4t}=\gamma_{st}=\gamma_{ts}=(\gamma_s\times\gamma_{ss})_s=(e_4\times e_{4s})_s=e_4\times e_{4ss},
\end{eqnarray*}%
which is exactly Eq.(\ref{flows6}). Thus we have showed the following
\begin{Theorem}
Eq.(\ref{2}) in Im$(\mathbb O)\cong\mathbb R^7$ is  equivalent to Eq.(\ref{flows6}) of Schr\"odinger flows from $\mathbb{R}^1$ into $\mathbb S^6$.
\end{Theorem}

We shall transform Eq.(\ref{flows6}), and hence Eq.(\ref{2}), to a nonlinear
Schr\"odinger-type system, like the Da Rios-NLS correspondence. In fact, from Eq.(\ref{2}), the fact: $e_{1s}=-\overline{\Arph}_1e_4+\Arph_2e_2$ and  Table \ref{table2}, we have that
\begin{eqnarray*}
\begin{split}
e_{4t}&=\gamma_{st}=\gamma_{ts}=(-\sqrt{-1}\Arph_1e_1+\sqrt{-1}\,\,\overline{\Arph}_1\overline{e}_1)_s\\
&=-\sqrt{-1}\,\Arph_{1s}e_1+\sqrt{-1}\,\overline{\Arph}_{1s}\overline{e}_1-\sqrt{-1}\,\Arph_1\Arph_2e_2+\sqrt{-1}\,\overline{\Arph}_1\overline{\Arph}_2\overline{e}_2.
\end{split}
\end{eqnarray*}
Hence, the complexified $G_2$-frame by $(e_4 \,\,\,f\,\,\,\overline{f}\,\,)^T$ admits
 \begin{eqnarray*}
\begin{array}{c@{\hspace{-5pt}}l}
\left(\begin{array}{ccc;{2pt/2pt}ccc}
 e_4\\
 f\\
 \overline{f}
\end{array}\right)_t
=\left(
\begin{array}{@{}ccc@{}}
0 & \omega &\overline{\omega} \\
-\overline{\omega}^T & \kappa & [\omega] \\
-\omega^T & \overline{[\omega]} & \overline{\kappa}
\end{array}
\right)
\end{array}
\left(\begin{array}{ccc;{2pt/2pt}ccc}
 e_4\\
 f\\
 \overline{f}
\end{array}\right),
\end{eqnarray*}
where, $\omega=(\omega^1\,\,\,\, \omega^2\,\,\,\,\omega^3)=(-\sqrt{-1}\,\Arph_{1s}\,\,\,\,-\sqrt{-1}\,\Arph_1\Arph_2\,\,\,\,0 ),\, R_1+R_2+R_3=0,$
\begin{eqnarray*}
\kappa=
\left(\begin{array}{@{}ccc@{}}
\sqrt{-1}R_1 & a_1 & a_2 \\
-\overline{a}_1 & \sqrt{-1}R_2 & a_3 \\
-\overline{a}_2 & -\overline{a}_3 & \sqrt{-1}R_3
\end{array}
\right),
\end{eqnarray*}%
\begin{eqnarray*}
\begin{split}
[\omega]=\frac{1}{\sqrt{2}}
\left(
\begin{array}{@{}ccc@{}}
0 & 0 & \Arph_1\Arph_2 \\
0 & 0 & -\Arph_{1s} \\
-\Arph_1\Arph_2 & \Arph_{1s} & 0
\end{array}
\right),
\end{split}
\end{eqnarray*}
and $R_i\in\mathbb{R}, a_i\in\mathbb{C},\,\,i\in\{1,2,3\}$ are functions with respect to $s$ and $t$, which will be determined later.

One the other hand, Eq.(\ref{fuhaibasis}) can be rewritten as
\begin{eqnarray*}
\begin{array}{c@{\hspace{-5pt}}l}
\left(\begin{array}{ccc}
 e_4\\
 f\\
 \overline{f}
\end{array}\right)_s
=\left(
\begin{array}{@{}ccc@{}}
0 & g & \overline{g} \\
-\overline{g}^T & M & [G] \\
-g^T & \overline{[G]} & \overline{M}
\end{array}
\right)
\end{array}
\left(\begin{array}{ccc;{2pt/2pt}ccc}
 e_4\\
 f\\
 \overline{f}
\end{array}\right),
\end{eqnarray*}
where
\begin{eqnarray*}
g=\left(\begin{array}{cc}
 \Arph_1\\
 0\\
 0
\end{array}\right)^T,\,\,\,\,
[G]=\left(\begin{array}{cccc}
0 & 0 & 0 \\
0 & 0 & -\frac{\sqrt{-1}}{\sqrt{2}}\Arph_1 \\
0 & \frac{\sqrt{-1}}{\sqrt{2}}\Arph_1 & 0
\end{array}\right),\,\,\,\,
M=\left(\begin{array}{@{}ccc@{}}
0 & \Arph_2 & 0 \\
-\overline{\Arph}_2 & 0 & \Arph_3 \\
0 & -\overline{\Arph}_3 & 0
\end{array}\right).
\end{eqnarray*}
From the integrability condition:
$\left(\begin{array}{ccc;{2pt/2pt}ccc}
e_4\\
f\\
\overline{f}
\end{array}\right)_{ts}=\left(\begin{array}{ccc;{2pt/2pt}ccc}
e_4\\
f\\
\overline{f}
\end{array}\right)_{st}$, we have:
\begin{subequations}
\begin{align}
g_t&=\omega_s+\omega M+\overline{\omega}\,\overline{[G]}-g\kappa-\overline{g}\,\overline{[\omega]},\label{gt1}\\
M_t&=\kappa_s-\overline{\omega}^Tg+\kappa M+[\omega]\,\overline{[G]}+\overline{g}^T\omega-M\kappa-[G]\,\overline{[\omega]},\label{gt2}\\
{[G]}_t&={[\omega]}_s-\overline{\omega}^T\,\overline{g}+\kappa [G]
+[\omega]\,\overline{M}+\overline{g}^T\,\overline{\omega}-M[\omega]-[G]\,\overline{\kappa},\label{gt3}
\end{align}
\end{subequations}
From (\ref{gt1}), (\ref{gt2}) and (\ref{gt3})\,\, we have
\begin{eqnarray}
\left\{\begin{array}{clll}\label{psi1}
\Arph_{1t}&=&-\sqrt{-1}\Arph_{1ss}+\sqrt{-1}\Arph_1\Arph_2\overline{\Arph}_2-\sqrt{-1}R_1\Arph_1,\\
a_1&=&-\frac{\sqrt{-1}}{\Arph_1}(2\Arph_{1s}\Arph_2+\Arph_1\Arph_{2s}),\\
a_2&=&-\frac{1}{\Arph_1}(\sqrt{-1}\Arph_{1}\Arph_2\Arph_3+\sqrt{2}\,\,\overline{\Arph}_{1}\overline{\Arph}_1\overline{\Arph}_2),
\end{array}\right.
\end{eqnarray}
and
\begin{eqnarray}
\left\{\begin{array}{clll}\label{psi2}
\Arph_{2t}&=&a_{1s}-\frac{3}{2}\sqrt{-1}\Arph_1\overline{\Arph}_1\Arph_2+\sqrt{-1}\Arph_2(R_1-R_2)-a_2\overline{\Arph}_3,\\
\Arph_{3t}&=&a_{3s}+\sqrt{-1}\Arph_3(R_2-R_3)+\overline{\Arph}_2\,a_2,\\
\sqrt{-1}R_{1s}&=&\sqrt{-1}(\Arph_1\overline{\Arph}_1)_s-\Arph_2\overline{a}_1+\overline{\Arph}_2a_1,\\
\sqrt{-1}R_{2s}&=&-\frac{\sqrt{-1}}{2}(\Arph_1\overline{\Arph}_1)_s-\overline{\Arph}_2a_1+\Arph_2\overline{a}_1
+\overline{\Arph}_3a_3-\Arph_3\overline{a}_3,\\
\sqrt{-1}R_{3s}&=&-\frac{\sqrt{-1}}{2}(\Arph_1\overline{\Arph}_1)_s-\overline{\Arph}_3a_3+\Arph_3\overline{a}_3,\\
\Arph_2a_3&=&a_{2s}+a_1\Arph_3.
\end{array}\right.
\end{eqnarray}
One notes that $R_1+R_2+R_3=0$ in (\ref{psi2}) which is just compatible with the requirement. Furthermore,
from (\ref{psi1})\,\,(\ref{psi2}), we obtain
\begin{eqnarray*}
\left\{\begin{array}{clll}
\Arph_{1t}&=&-\sqrt{-1}\Arph_{1ss}+\sqrt{-1}\Arph_1\Arph_2\overline{\Arph}_2-\sqrt{-1}R_1\Arph_1,\\
\Arph_{2t}&=&a_{1s}-\frac{3}{2}\sqrt{-1}\Arph_1\overline{\Arph}_1\Arph_2+\sqrt{-1}\Arph_2(R_1-R_2)-a_2\overline{\Arph}_3,\\
\Arph_{3t}&=&a_{3s}+\sqrt{-1}\Arph_3(R_1+2R_2)+\overline{\Arph}_2\,a_2,\\
\end{array}\right.
\end{eqnarray*}
where
\begin{subequations}
	\begin{align}
a_1=&-\frac{\sqrt{-1}}{\Arph_1}(2\Arph_{1s}\Arph_2+\Arph_1\Arph_{2s}),\label{a11}\\
a_2=&-\frac{1}{\Arph_1}(\sqrt{-1}\Arph_{1}\Arph_2\Arph_3+\sqrt{2}\,\,\overline{\Arph}_{1}\overline{\Arph}_1\overline{\Arph}_2),\label{a22}\\
a_3=&-2\sqrt{-1}\Arph_3(\ln \Arph_1\Arph_2)_s-\sqrt{-1}\Arph_{3s}-\frac{\sqrt{2}}{\Arph_2}\,(\frac{\overline{\Arph}_1\overline{\Arph}_1\overline{\Arph}_2}{\Arph_1})_s,\label{a33}\\
\sqrt{-1}R_{1s}=&\sqrt{-1}(\Arph_1\overline{\Arph}_1)_s-\sqrt{-1}(\Arph_2\overline{\Arph}_2)_s-2\sqrt{-1}\Arph_2\overline{\Arph}_2(\ln \Arph_1\overline{\Arph}_1)_s, \label{a44}\\
\sqrt{-1}R_{2s}=&-\frac{\sqrt{-1}}{2}(\Arph_1\overline{\Arph}_1)_s+\sqrt{-1}(\Arph_2\overline{\Arph}_2)_s
-\sqrt{-1}(\Arph_3\overline{\Arph}_3)_s+2\sqrt{-1}\Arph_2\overline{\Arph}_2(\ln \Arph_1\overline{\Arph}_1)_s\nonumber\\
&-2\sqrt{-1}\Arph_3\overline{\Arph}_3(\ln \Arph_1\overline{\Arph}_1\Arph_2\overline{\Arph}_2)_s+\sqrt{2}\,[\,\frac{\Arph_3}{\overline{\Arph}_2}
(\frac{\Arph_1\Arph_1\Arph_2}{\overline{\Arph}_1})_s
-\frac{\overline{\Arph}_3}{\Arph_2}(\frac{\overline{\Arph}_1\overline{\Arph}_1\overline{\Arph}_2}{\Arph_1})_s]\label{a55}.
    \end{align}
\end{subequations}
The equations (\ref{a44}) and (\ref{a55}) imply that $R_1$ and $R_2$  are of the forms as follows
\begin{eqnarray*}
	\begin{split}
		\begin{array}{rll}
			R_1=&|\Arph_1|^2-|\Arph_2|^2-2\int_0^s |\Arph_2|^2(\ln |\Arph_1|^2)_sd\widetilde{s}+R_{10}(t),\\
			R_2=&\int_0^s \{2|\Arph_2|^2(\ln |\Arph_1|^2)_s-2|\Arph_3|^2(\ln (|\Arph_1|^2|\Arph_2|^2))_s-\sqrt{2}\sqrt{-1}[\frac{\Arph_3}{\overline{\Arph}_2}(\frac{\Arph^2_1\Arph_2}
{\overline{\Arph}_1})_s-\frac{\overline{\Arph}_3}{\Arph_2}(\frac{\overline{\Arph}^2_1\overline{\Arph}_2}{\Arph_1})_s]\}d\widetilde{s}\\
			&-\frac{1}{2}|\Arph_1|^2+|\Arph_2|^2-|\Arph_3|^2+R_{20}(t),
		\end{array}
	\end{split}
\end{eqnarray*}
where $R_{10}$ and $R_{20}$ depend only on $t$.
Now, by the following transformations:
\begin{eqnarray*}
\begin{split}
\Arph_1& \mapsto \sqrt{2}\,\Arph_1\exp (-\sqrt{-1}\int_0^t R_{10}d\widetilde{t}),\\
\Arph_2& \mapsto \Arph_2\exp (\sqrt{-1}\int_0^t (R_{10}-R_{20})d\widetilde{t}),\\
\Arph_3& \mapsto \Arph_3\exp (\sqrt{-1}\int_0^t (R_{10}+2R_{20})d\widetilde{t}),
\end{split}
\end{eqnarray*}
we arrive at
\begin{Theorem}\label{them2} If $\kappa_2\neq 0$, then the equation of Schr\"odinger flows to $\mathbb S^6 \subset \mathbb{R}^7$ is equivalent to
\begin{eqnarray}\label{NLST}
\left\{\begin{array}{lll}
\sqrt{-1}\Arph_{1t}&=&\Arph_{1ss}+2\Arph_1|\Arph_1|^2-2\Arph_1|\Arph_2|^2-2\Arph_1\int_0^s|\Arph_2|^2(\ln|\Arph_1|^2)_sd\widetilde{s},\\
\sqrt{-1}\Arph_{2t}&=&\Arph_{2ss}+2\Arph_2|\Arph_2|^2+2(\Arph_2(\ln \Arph_1)_s)_s-2\Arph_2|\Arph_3|^2+2\sqrt{-1}\,\frac{\overline{\Arph}_1^2\,\overline{\Arph}_2\overline{\Arph}_3}{\Arph_1}\\
&&+2\Arph_2\int_0^s\{2|\Arph_2|^2(\ln|\Arph_1|^2)_s-|\Arph_3|^2(\ln(|\Arph_1|^2|\Arph_2|^2))_s
+2\text{Im}\,[\,\frac{\Arph_3}{\overline{\Arph}_2}(\frac{\Arph_1^2\Arph_2}{\overline{\Arph}_1})_s\,]\}d\widetilde{s},\\
\sqrt{-1}\Arph_{3t}&=&\Arph_{3ss}+2\Arph_3|\Arph_3|^2+2[\,\Arph_3(\ln (\Arph_1\Arph_2))_s-\frac{\sqrt{-1}}{\Arph_2}(\frac{\overline{\Arph}_1^2\overline{\Arph}_2}{\Arph_1})_s\,]_s-2\sqrt{-1}\frac{\overline{\Arph}_1^2\overline{\Arph}_2^2}{\Arph_1}\\
&&+2\Arph_3\int_0^s\{2|\Arph_3|^2(\ln|\Arph_1|^2|\Arph_2|^2)_s-|\Arph_2|^2(\ln|\Arph_1|^2)_s
-4\,\text{Im}\,[\,\frac{\Arph_3}{\overline{\Arph}_2}(\frac{\Arph_1^2\Arph_2}{\overline{\Arph}_1})_s\,]\}d\widetilde{s},
\end{array}\right.
\end{eqnarray}
which is a nonlinear Schr\"odinger-type system (NLSS) in three unknown complex functions.
\end{Theorem}

\begin{Remark}
From the proof of Theorem \ref{them2}, one sees that it is by an application of the $G_2$ structure displayed in Theorem \ref{th-ab} that
we obtain the nonlinear Schr\"odinger-type system (NLSS) (\ref{NLST}). Moreover, one also notes that the same sign in the terms $e_1$ and $e_2$ in the equation $e_{4t}=e_4\times e_{4ss}$  produces an integration term in $a_1$ (see (\ref{a11})) and hence in the first equation in (\ref{NLST}). The additional integration terms in $a_2$ and $a_3$ (see (\ref{a22}) and (\ref{a33})) (i.e. the second and third equations in (\ref{NLST})) are in fact produced by the $G_2$-structure itself. Moreover,
if $\kappa_2=0$, i.e. $\Arph_2=0,$  Eq.(\ref{psi1}) and (\ref{psi2}) are clearly reduced to
\begin{eqnarray*}
\left\{\begin{array}{clll}
\sqrt{-1}\Arph_{1t}&=&\Arph_{1ss}+2\Arph_1|\Arph_1|^2,\\
\sqrt{-1}\Arph_{3t}&=&\Arph_{3ss}+2\Arph_3|\Arph_3|^2,
\end{array}\right.
\end{eqnarray*}
by choosing $a_3=-\sqrt{-1}\varphi_{3s}$, and hence Eq.(\ref{NLST}) is reduced to the focusing nonlinear Schr\"odinger equation (NLS).
\end{Remark}

Conversely, if $\Arph_1(t,s)$, $\Arph_2(t,s)$ and $\Arph_3(t,s)$ are a solution to Eq.(\ref{NLST}), one may
verify that there is a corresponding family of curves $\gamma(t,s)$ in $\mathbb R^7$ satisfying Eq.(\ref{2}). The details are omitted here.
This produces Eq.(\ref{2})-NLSS correspondence, that is, if a curves $\gamma(t,s)$ evolves according to Eq.(\ref{2}), then the associated complex functions $\varphi_1, \varphi_2$ and $\varphi_3$ given by (\ref{Arph1}) evolve according to Eq.(\ref{NLST}). This generalizes the Da Rios-NLS correspondence mentioned in the introduction, since when $\varphi_2=0$ the correspondence is reduced to the Da Rios-NLS correspondence.
One realizes that though Kirchhoff's almost complex structure on $\mathbb S^6$ is not integrable, the  Eq.(\ref{2})-NLSS correspondence is still valid.

It seems that it is the non-integrability of Kirchhoff's almost complex structure on $\mathbb S^6$ that creates the integration terms in NLSS (\ref{NLST}), as there is no an integration term involved in the NLS equation. If there is an almost complex structure $J$ on $\mathbb S^6$
for which the corresponding nonlinear Schr\"odinger-type system is without integration terms, we believe such an almost complex structure is  integrable. Does there exist such an almost complex structure on $\mathbb S^6$?

To complete this section, we come to find an almost complex structure structure on $\mathbb S^6 $, which is different from Kirchhoff's almost complex structure, such that the integration terms in the corresponding nonlinear Schr\"odinger-like system are reduced greatly.

For $A\in O(7)/G_2$, Calabi constructed in \cite{Calabi} an almost complex structure $J^A$ on $\mathbb S^6$, which is still compatible to the
standard round metric on $\mathbb S^6$. The almost complex structure $J^A$ reads
\begin{eqnarray*}
J^A: T_u \mathbb S^6 \rightarrow T_u \mathbb S^6, \,\,
X\mapsto J^A(X)=A^{-1}((Au)\times (AX)),~ X\in T_u\mathbb S^6.
\end{eqnarray*}
One sees that only when $A\in G_2$, $J^A$ coincides with Kirchhoff's almost complex structure $J$ on $\mathbb S^6$. In fact, in this case
for $u\in \mathbb S^6$,
$$J_u^A(X)=A^{-1}((Au)\times (AX))=A^{-1}(Au)\times A^{-1}(AX)=u\times X=J_{u}(X),~ \forall X\in  T_u\mathbb S^6.$$
With the almost complex structure $J^A$, one may verify straightforwardly that the equation of Schr\"odinger flows from $\mathbb R$ to $\mathbb S^6\hookrightarrow\mathbb R^{7}$ is
\begin{eqnarray}\label{flow12}
	u_t=A^{-1}(Au\times(A(u_{ss}))),~~ u: \mathbb R\times \mathbb R\to \mathbb S^6\hookrightarrow\mathbb R^{7}.\nonumber
\end{eqnarray}
Correspondingly, Eq.(\ref{2}) is modified to be
\begin{eqnarray}\label{e4A1}
\gamma_{4t}=A^{-1}(A\gamma_s\times (A\gamma_{ss}),
\end{eqnarray}
where $\gamma=\gamma(t,s)\in\mathbb R^7$.

Now by choosing
\begin{eqnarray*}
A=\left(
\begin{array}{@{}ccc|cc|cc@{}}
1 & 0 & 0 & 0 & 0 & 0 & 0 \\
0 & 1 & 0 & 0 & 0 & 0& 0 \\
0 & 0 & 1 & 0 & 0 & 0 & 0 \\
\hline
0 & 0 & 0 & 0 & 1 & 0 & 0 \\
0 & 0 & 0 & -1 & 0 & 0 & 0 \\
\hline
0 & 0 & 0 & 0 & 0 & 1 & 0 \\
0 & 0 & 0 & 0 & 0 & 0 & 1 \\
\end{array}
\right)\in O(7)/G_2
\end{eqnarray*}
and with a direct but long computation similar to that in the proof of Theorem \ref{them2}, we have
\begin{Proposition}
Eq.(\ref{e4A1}) is equivalent to the following nonlinear Schr\"odinger-like system in three unknown complex-functions:
\begin{eqnarray}\label{1111}
\left\{\begin{array}{lll}
\sqrt{-1}\Arph_{1t}&=&\Arph_{1ss}+2\Arph_1|\Arph_1|^2+2\Arph_1|\Arph_2|^2.\\
\sqrt{-1}\Arph_{2t}&=&-\Arph_{2ss}-2\Arph_2|\Arph_2|^2-6\Arph_2|\Arph_1|^2+2\Arph_2|\Arph_3|^2-2\sqrt{-1}\,\frac{\overline{\Arph}_1^2\,\overline{\Arph}_2\overline{\Arph}_3}{\Arph_1}\\
&&+2\Arph_2\int_0^s\{|\Arph_3|^2(\ln(|\Arph_2|^2))_s
-2\text{Im}\,[\,\frac{\Arph_3}{\overline{\Arph}_2}(\frac{\Arph_1^2\Arph_2}{\overline{\Arph}_1})_s\,]\}d\widetilde{s},\\
\sqrt{-1}\Arph_{3t}&=&-\Arph_{3ss}-2\Arph_3|\Arph_3|^2-2[\,\Arph_3(\ln \Arph_2)_s-\frac{\sqrt{-1}}{\Arph_2}(\frac{\overline{\Arph}_1^2\overline{\Arph}_2}{\Arph_1})_s\,]_s+2\sqrt{-1}\frac{\overline{\Arph}_1^2\overline{\Arph}_2^2}{\Arph_1}\\
&&-4\Arph_3\int_0^s\{2|\Arph_3|^2(\ln|\Arph_2|^2)_s
-2\,\text{Im}\,[\,\frac{\Arph_3}{\overline{\Arph}_2}(\frac{\Arph_1^2\Arph_2}{\overline{\Arph}_1})_s\,]\}d\widetilde{s}.
\end{array}\right.
\end{eqnarray}
\end{Proposition}

We would point out that it is due to the choice of $A$ indicated above that the last term in (\ref{a44}) is cancelled and hence the integration terms in the corresponding nonlinear Schr\"odinger-like system (\ref{1111}) are reduced.

\section *{\S 5 Discussions and remarks}

In order to exploit further geometric properties of $\Arph_1,\Arph_2$
and $\Arph_3$, we introduce a surface $\Sigma$ $ \hookrightarrow \mathbb{R}^7$= Im $ \mathbb O $ swept by
$\gamma(t,s).$  Let $\{E_1=\gamma_s=I_4,\, E_2=\frac{\gamma_t}{k_1}=-I_5,\, E_3=I_1,\, E_4=I_2,\,E_5=I_6,\, E_6=I_3,\,E_7=I_7 \}$ be an orthonormal $G_2$-frame associated to $\Sigma$ in Im($\mathbb O)\cong\mathbb R^7$ such that $\{E_1, E_2\}$ spans
$T_{(t,s)}\Sigma$, $\{E_3,\,E_4,\,E_5,\,E_6,\,E_7\}$ spans $(T_{(t,s)}\Sigma)^\perp$, and $\{\omega^1,\cdots,\omega^7\}$ is its co-frame.
The structure equations of $\Sigma$ are given by
\begin{eqnarray}\label{structure1}
\left\{\begin{array}{lll}
d\gamma=\omega^i\,E_i,\\
dE_i=\omega^k_i\,E_k+\omega_i^\alpha\,E_{\alpha},\\
dE_{\alpha}=\omega^k_{\alpha}\,E_k+\omega_{\alpha}^\beta\,E_{\beta},\\
d\omega^i=-\omega_j^i\wedge\omega^j=0,\, \omega^i_j+\omega_i^j,\\
d\omega^\alpha=-\omega_i^\alpha\wedge\omega^i,
\end{array}\right.
\end{eqnarray}
where $i,j,k=1,2$ and $\alpha,\beta,\eta=3,4,5,6,7$. Restricting to $\Sigma$, we have $\omega^\alpha=0$. Hence, from
\begin{eqnarray*}
0=d\omega^\alpha=-\omega^\alpha_i\wedge\omega^i,
\end{eqnarray*}
we have
\begin{eqnarray}\label{hij1}
\omega_i^\alpha=h_{ij}^\alpha\omega^j,\,~~h_{ij}^\alpha=h_{ji}^\alpha.
\end{eqnarray}
One knows that $\{h_{ij}^{\alpha}\}$ ($1\le i,j\le2, 3\le\alpha\le7$) are the coefficients of the second fundamental form of $\Sigma\hookrightarrow \mathbb R^7\cong {\rm Im}(\mathbb O)$.

\begin{Lemma} \label{lemma1}
The non-zero coefficients of the second fundamental form $(h^\alpha_{ij})$ are
\begin{eqnarray}\label{second1}
\begin{split}
h^3_{11}=&\sqrt{2}|\Arph_1|\,,h^3_{12}=h_{21}^3=\frac{\sqrt{-1}}{\sqrt{2}}(\ln \frac{\overline{\Arph}_1}{\Arph_1})_s,\\
h_{22}^3=&-\frac{\Arph_1\overline{\Arph}_{1ss}+\overline{\Arph}_1\Arph_{1ss}}{2\sqrt{2}\,|\Arph_1|^3}+\frac{|\Arph_2|^2}{\sqrt{2}\,|\Arph_1|},\\
h^4_{22}=&-\frac{2(\ln |\Arph_1|)_s|\Arph_2|+(|\Arph_2|)_s}{\sqrt{2}\,|\Arph_1|},\\
h_{22}^5=&-\frac{\sqrt{-1}\,[2(\ln \frac{\overline{\Arph}_1}{\Arph_1})_s\,|\Arph_2|^2+\Arph_2\overline{\Arph}_{2s}-\overline{\Arph}_2\Arph_{2s}]}{2\sqrt{2}\,|\Arph_1||\Arph_2|},\\
h^6_{22}=&\frac{\Arph_1^3\Arph_2^2\Arph_3+\overline{\Arph}_1^3\overline{\Arph}_2^2\overline{\Arph}_3}{2\sqrt{2}\,|\Arph_1|^4|\Arph_2|},\,h_{12}^5=h_{21}^5=-|\Arph_2|,\\
h^7_{22}=&\frac{\sqrt{-1}(\overline{\Arph}_1^3\overline{\Arph}_2^2\overline{\Arph}_3-\Arph_1^3\Arph_2^2\Arph_3)}{2\sqrt{2}\,|\Arph_1|^4|\Arph_2|}-\frac{3|\Arph_2|}{2}.
\end{split}
\end{eqnarray}
\end{Lemma}

{\bf Proof}:
The first fundamental form of $\gamma(t,s)$ is
\begin{eqnarray*}
\begin{split}
ds_{M}^2&=\langle \frac{\partial}{\partial s},\frac{\partial}{\partial s} \rangle ds^2+\langle \frac{\partial}{\partial t},\frac{\partial}{\partial t} \rangle dt^2=\langle \gamma_s,\gamma_s \rangle ds^2+\langle  \gamma_t,\gamma_t \rangle dt^2\\
&=ds^2+2|\Arph_1|^2dt^2=(\omega^1)^2+(\omega^2)^2,
\end{split}
\end{eqnarray*}
where $\omega^1=ds$ and $\omega^2=\sqrt{2}|\Arph_1|dt$. By Eq.(\ref{basis}) and Theorem \ref{thembasis}, a straightforward
computation shows that
\begin{eqnarray}\label{j1}
\begin{split}
dE_1=\frac{(|\Arph_1|)_s}{|\Arph_1|}\omega^2E_2+[\sqrt{2}|\Arph_1|\omega^1+\frac{\sqrt{-1}}{\sqrt{2}}(\ln\frac{\overline{\Arph}_1}{\Arph_1})_s\omega^2\,]E_3-|\Arph_2|\omega^2E_5,\,\,\,\,\,\,\,\,\,\,\,\,\,\,\,\,\,\,\,\,\,\,\,\,\,\,\,\,\,\,\,\,\,\,\,
\end{split}
\end{eqnarray}
and
\begin{eqnarray}\label{j2}
\begin{split}
dE_2=&-dI_5=d(-\sqrt{-1}\,\overline{r}e_1+\sqrt{-1}\,r\overline{e}_1)\\
=&-\frac{(|\Arph_1|)_s}{|\Arph_1|}\omega^2E_1+[\frac{\sqrt{-1}}{\sqrt{2}}(\ln \frac{\overline{\Arph}_1}{\Arph_1})_s \omega^1+(-\frac{\Arph_1\overline{\Arph}_{1ss}+\overline{\Arph}_1\Arph_{1ss}}{2\sqrt{2}\,|\Arph_1|^3}+\frac{|\Arph_2|^2}{\sqrt{2}\,|\Arph_1|})\omega^2]E_3\\
&-\frac{2(\ln |\Arph_1|)_s|\Arph_2|+(|\Arph_2|)_s}{\sqrt{2}\,|\Arph_1|}\omega^2E_4+\frac{\Arph_1^3\Arph_2^2\Arph_3+\overline{\Arph}_1^3\overline{\Arph}_2^2\overline{\Arph}_3}{2\sqrt{2}\,|\Arph_1|^4|\Arph_2|}\omega^2E_6\\
&-[\,|\Arph_2|\omega^1+\frac{\sqrt{-1}\,(2(\ln \frac{\overline{\Arph}_1}{\Arph_1})_s\,|\Arph_2|^2+\Arph_2\overline{\Arph}_{2s}-\overline{\Arph}_2\Arph_{2s})}{2\sqrt{2}\,|\Arph_1||\Arph_2|}\omega^2\,]E_5\\
&+[\frac{\sqrt{-1}(\overline{\Arph}_1^3\overline{\Arph}_2^2\overline{\Arph}_3-\Arph_1^3\Arph_2^2\Arph_3)}{2\sqrt{2}\,|\Arph_1|^4|\Arph_2|}-\frac{3|\Arph_2|}{2}\omega^2]E_7.
\end{split}
\end{eqnarray}
Since
\begin{eqnarray}\label{j3}
 dE_i=\sum\limits_{k=1}^2\omega^k_i\,E_k+\sum\limits_{\alpha=3}^7(h_{i1}^\alpha\omega^1+h_{i2}^\alpha\omega^2)\,E_{\alpha},
\end{eqnarray}
based on Eqs.(\ref{structure1}) and (\ref{hij1}), we have (\ref{second1}) from
Eqs.(\ref{j1}), (\ref{j2}) and (\ref{j3}).
\qed

When $\Arph_2=0$, $h_{ij}^\alpha=0$ for $ \alpha= 4, 5, 6, 7$ indicates that $\Sigma$ is located completely in an associative 3-dimensional space $ \mathbb{V}^3=\text{span}_{\mathbb{R}}\{I_4, I_1, I_5\} $. Therefore, Eq.(\ref{NLST}) is reduced to the nonlinear Schr\"odinger equation (NLS)
\begin{eqnarray*}\label{NLS2}
\sqrt{-1}\Arph_{1t}=\Arph_{1ss}+2\Arph_1|\Arph_1|^2.
\end{eqnarray*}
In this case, it can be inferred that the almost complex structure $J$ on $\mathbb S^6$ returns to the complex (integrable) structure $ J $ on $\mathbb S^2$.
We note that a 3-dimensional vector space $ \mathbb{V} $ in Im($\mathbb O)=\mathbb{R}^7$ is called associative if
$ \mathbb{V}=\text{span}_{\mathbb{R}}\{u, v, u\times v\} $. Associative 3-planes or co-associative 4-planes in Im($\mathbb O)$ are well-studied in $G_2$-geometry. This gives the geometric characterization of the surface $\Sigma$ corresponding to the case that $\varphi_2=0$.

In order to characterize the geometric properties corresponding to $\Arph_3=0$ and $\varphi_2\not=0$, we set $ \widetilde{E}_i=E_i $ except that $ \widetilde{E}_4=\cos\theta E_4-\sin\theta E_7$ and $\widetilde{E}_7=\sin\theta E_4+\cos\theta E_7, $ where $\theta$ will be determined later. Then $\{\widetilde{E}_1,\,\widetilde{E}_2,\,\widetilde{E}_3,\,\widetilde{E}_4,\,\widetilde{E}_5,\,\widetilde{E}_6,\,\widetilde{E}_7\,\}$ also consists of an orthonormal $SO(7)$-frame associated to $\Sigma$ in Im($\mathbb O)\cong \mathbb R^7$. Let $\{\widetilde{\omega}^1,\cdots,\widetilde{\omega}^7\}$ be its dual-frame. Since $\widetilde{E}_1=E_1,\,\widetilde{E}_2=E_2,$ we have $\widetilde{\omega}^1=\omega^1,\,\widetilde{\omega}^2=\omega^2.$

\begin{Lemma} \label{lemma2} By choosing $\theta=\arccos\frac{h_{22}^4}{\sqrt{(h_{22}^4)^2+\frac{9}{4}|\Arph_2|^2}},$ the coefficients of the second fundamental form $ \widetilde{h}^\alpha_{ij} $ corresponding to the orthonormal frame $\{\widetilde{E}_1,\cdots, \widetilde{E}_7\}$ are
\begin{eqnarray*}
\begin{split}
\widetilde{h}_{ij}^3=&h_{ij}^3,\,\,\widetilde{h}_{ij}^5=h_{ij}^5,\,\,\widetilde{h}_{ij}^6=h_{ij}^6,\\
\widetilde{h}_{11}^4=&\widetilde{h}_{12}^4=\widetilde{h}_{21}^4=0,\\ \widetilde{h}_{22}^4=&\sqrt{(h_{22}^4)^2+\frac{9}{4}|\Arph_2|^2}-\frac{3|\Arph_2|(h_{22}^7+\frac{3}{2}|\Arph_2|)}{2\sqrt{(h_{22}^4)^2+\frac{9}{4}|\Arph_2|^2}},\\
\widetilde{h}_{11}^7=&\widetilde{h}_{12}^7=\widetilde{h}_{21}^7=0,\widetilde{h}_{22}^7=\frac{h_{22}^4(h_{22}^7+\frac{3}{2}|\Arph_2|)}{\sqrt{(h_{22}^4)^2+\frac{9}{4}|\Arph_2|^2}}.
\end{split}
\end{eqnarray*}
\end{Lemma}

{\bf Proof}: By a direct computation.\qed

\begin{Proposition}\label{prop1}
Let $\Arph_3=0$ and $\Arph_2\neq0$. Then the coefficients of the second fundamental form are: ${\widetilde h}^{\alpha}_{ij}=0$ for $\alpha\in\{6,7\}$, namely, the normal bundle of the surface $ \Sigma $ is flat in directions  $E_i$ ($i=6,7$).
\end{Proposition}

{\bf Proof}:
From Lemma \ref{lemma1}, Lemma \ref{lemma2} and $\Arph_3=0$, the corresponding nonzero terms for $\alpha\in\{6,7\}$ now become
\begin{eqnarray*}
\begin{split}
\widetilde{h}_{11}^6=&\widetilde{h}_{12}^6=\widetilde{h}_{21}^6=\widetilde{h}_{11}^7=\widetilde{h}_{12}^7=\widetilde{h}_{21}^7=0,\,
\widetilde{h}^6_{22}=\frac{\Arph_1^3\Arph_2^2\Arph_3+\overline{\Arph}_1^3\overline{\Arph}_2^2\overline{\Arph}_3}{2\sqrt{2}\,|\Arph_1|^4|\Arph_2|}=0,\\
\widetilde{h}^7_{22}=&\frac{\sqrt{-1}(\overline{\Arph}_1^3\overline{\Arph}_2^2\overline{\Arph}_3-\Arph_1^3\Arph_2^2\Arph_3)}{2\sqrt{2}\,|\Arph_1|^4|\Arph_2|}\frac{h_{22}^4}{\sqrt{(h_{22}^4)^2+\frac{9}{4}|\Arph_2|^2}}=0.
\end{split}
\end{eqnarray*}
This completes the proof of Lemma \ref{prop1}.\qed

\medskip
The geometric properties of $\varphi_i$ ($i=1,2,3$) are summarized as follows. When $\Arph_2=0$, the surface $\Sigma$ (swept by
$\gamma(t,s)$) is located completely in a 3-manifold spanned by $\{E_1,E_2,E_3\}$ with the algebraic property of associativity;
when $\Arph_3=0$ and $\varphi_2\not=0$, the normal bundle of the surface $ \Sigma $ is flat in directions $E_i$ ($i=6,7$) and roughly speaking,
we can regarded $\Sigma$ as located in a 5-manifold spanned by
$\{\widetilde{E}_1,\,\widetilde{E}_2,\,\widetilde{E}_3,\,\widetilde{E}_4,\,\widetilde{E}_5\}$.

\medskip
We have furthered our understanding of almost complex structures on $\mathbb S^6 $ and the $G_2$-structure on Im$(\mathbb O) = \mathbb{R}^7 $ via Schr\"odinger flows. Many related problems remain to be clarified. For example, the 5-manifold spanned by $\{\widetilde{E}_1,\,\widetilde{E}_2,\,\widetilde{E}_3,\,\widetilde{E}_4,\,\widetilde{E}_5\}$ corresponding to $\Arph_3=0$ and $\varphi_2\not=0$ is
a new object in relation to the surface $\Sigma$ swept by $\gamma(t,s)$ in $\mathbb R^7$. How might one further characterize its geometric or algebraic properties? It is very interesting to note that like the usual binormal motion equation (\ref{1}) in $\mathbb R^3$, whether the $G_2$-binormal motion equation (\ref{2}) in $\mathbb R^7$ is integrable, and furthermore, whether it sits in a hierarchy with other new integrable systems and has relations to various known integrable systems.

\section * {\bf Acknowledgements}
The authors would like to express their deep thanks to the referee for his/her helpful comments.
This work is supported by the National Natural Science Foundation of
China (Grant No. 11271073).


\begin{thebibliography}{399}

\bibitem{Arms}  \label{Arms}
R.J. Arms and F.R. Hama, Localized-induction concept on a curved vortex and motion of elliptic vortex
ring, Phys. Fluids {\bf 8} (1965), 553-559.

\bibitem{Baez} \label{Baez}		
J.C. Baez, The octonions, Bull. Amer. Math. Soc. {\bf 39} (2002), 145-206.	

\bibitem{Borel} \label{Borel}
A. Borel  and J.P. Serre, Groupes de Lie et puissances r\'eduites de Steenrod, Amer. J. Math. {\bf 75} (1953), 409-448.

\bibitem{Brower1}  \label{Brower1}
B.C. Brower, D.A. Kessler, J. Koplik  and H. Levine, Geometrical models of interface evolution, Phys. Rev. A {\bf 29} (1984), 1335-1342.

\bibitem{Bryant} \label{Bryant}
R.L Bryant, Submanifolds and special structures on the octonions, J. Diff. Geom. {\bf 17} (1982), 185-232.

\bibitem{Bryant11}  \label{Bryant11}
R.L. Bryant, S.S. Chern's study
of almost-complex structure on the six-sphere, arXiv: 14053405.

\bibitem{Calabi} \label{Calabi}
E. Calabi, Construction and properties of some 6-dimensional almost complex
manifolds, Trans. A.M.S. {\bf 87} (1958), 407-438.

\bibitem{Gluck} \label{Gluck}
E. Calabi  and H. Gluck, What are the best almost complex structures on the 6-sphere, Proceedings of Symposia in Pure Mathematics {\bf 54} (1993), 99-106.

\bibitem{Rios} \label{Rios}
L.S. Da Rios, On the motion of an unbounded fluid with a vortex filament of any shape, Rend. Circ. Mat. Palermo. {\bf 22} (1906), 117-135.

\bibitem{dq} \label{dq}
Q. Ding, A note on the NLS and the Schr\"odinger flow of maps, Phys.
Lett. A {\bf 248} (1998), 49-56.

\bibitem{dinghe}
Q. Ding and Z.Z. He, The noncommutative KdV equation and its
para-K\"ahler structure, Sci. in China Math. {\bf 57} (2014)
1505-1516.

\bibitem{DW1}
Q. Ding and Y.D. Wang, Geometric KdV flows, motions of curves and
the third order system of the AKNS hierarchy, Inter. J. of Math.
{\bf 22} (2011) 1013-1029.

\bibitem{dw} \label{dw}
W.Y. Ding and Y.D. Wang, Schr\"odinger flows of maps into symplectic manifolds, Sci. China A {\bf 41} (1998), 746-755.

\bibitem{Harvey}  \label{Harvey}
R. Harvey and H.B. Lawson, Calibrated geometries, Acta Math. {\bf 148} (1982), 47-157.

\bibitem{Hasimoto} \label{Hasimoto}	
H. Hasimoto, A soliton on a vortex filament, J. Fluid. Mech. {\bf 51} (1972), 477-485.

\bibitem{Hashimoto} \label{Hashimoto}
H. Hashimoto and M. Ohashi, Orthogonal almost complex structures of hypersurfaces of purely imaginary octonions, Hokkaido Math. J. {\bf 39} (2010), 351-387.

\bibitem{LangerJS}  \label{LangerJS}
J.S. Langer, Instabilities and pattern formation in crystal growth, Rev. Mod. Phys. {\bf 52} (1980), 1-28.

\bibitem{Langer} \label{Langer}
J. Langer and R. Perline, Geometric realizations of Fordy-Kulish nonlinear Schr\"odinger sysyems, Pacific
J. Math. {\bf 195} (2000), 157-178.

\bibitem{LeBrun}  \label{LeBrun}
C. LeBrun, Orthogonal complex structures on $\mathbb S^6 $, Proc. A.M.S. {\bf 101} (1987), 136-138.

\bibitem{Meleshko}  \label{Meleshko}
V.V. Meleshko, A.A. Gourjii and T.S. Krasnopolskaya, Vortex ring: history and state of the art, J. of
Math. Sciences {\bf 187} (2012), 772-806.

\bibitem{Nahmod} \label{Nahmod}
A. Nahmod, A. Stefanov and K. Uhlenbeck, On Schr\"oinger maps,
Comm. Pure Appl. Math. {\bf 56} (2003),  114-151; Erratum,
On Schr\"odinger maps, Comm. Pure Appl. Math. {\bf 57} (2004), 833-839.

\bibitem{ohashi2011} \label{ohashi2011}
M. Ohashi, On $G_2$-invariants of curves of purely imaginary octonions, Recent Progress in Differential Geometry and Its Related Fields, 2012, 25-40.

\bibitem{ohashi20132} \label{ohashi20132}
M. Ohashi, $G_2$-congruence theorem for curves in purely imaginary octonions and its application, Geometriae Dedicata {\bf 163} (2013), 1-17.

\bibitem{ohashi2015}  \label{ohashi2015}
M. Ohashi, A method of derermining the SO(7)-invariants for curves in Im($\mathbb O$) by their $G_2$-invariants, Current Developments in Differential Geometry and its Related Fields, 2016, 201-213.

\bibitem{Saffman}  \label{Saffman}
P.G. Saffman and G. Taylor, The penetration of a fluid into a porous medium or Hale-Shaw cell
containing a more visous, Proc. R. Soc. London A {\bf 245} (1958), 312-329.

\bibitem{TeUh}
C.L. Terng and K. Uhlenbeck, Schr\"odinger flows on Grassmannians,
AMS/IP Studies in Advanced Mathematics, Vol. 36 (2006) 235--256.


\end{thebibliography}
\end{document}